\documentclass[11pt]{article}
\usepackage{latexsym, amscd, amsfonts, eucal, mathrsfs, amsmath, amssymb, amsthm, xypic,xr, stmaryrd, color, enumerate, tikz}
\usepackage{mathabx}
\usepackage{appendix}
\usepackage[all]{xy}
\usepackage{hyperref}
\usepackage{fullpage}
\newtheorem{theorem}{Theorem}[subsection]

\newtheorem{proposition}[theorem]{Proposition}
\newtheorem{corollary}[theorem]{Corollary}

\theoremstyle{definition}
\newtheorem{definition}[theorem]{Definition}
\newtheorem{construction}[theorem]{Construction}

\newtheorem{remark}[theorem]{Remark}

\newtheorem*{convention}{Convention}

\usepackage{cleveref}
\usepackage{scrextend}

\DeclareMathOperator*{\holim}{holim}

\begin{document}
\title{$C_2$-equivariant Homology Operations: Results and Formulas}
\author{Dylan Wilson}
\maketitle
\begin{abstract} 
	In this note we state corrected and expanded versions of
	our previous results on
	power operations for $C_2$-equivariant Bredon homology
	with coefficients in the constant Mackey functor on $\mathbb{F}_2$.
	In particular, we give a version of the Adem relations. The proofs
	rely on certain results in
	equivariant higher algebra which we will supply in a longer
	version of this paper.
\end{abstract}
\newpage
\tableofcontents
\newpage

\section*{Introduction}\label{sec:intro}
\addcontentsline{toc}{section}{\nameref{sec:intro}}

If $X$ is a spectrum then the smash power $X^{\wedge 2}$ carries
a $\Sigma_2$-action. Equipping $X$ with a map
	\[
	X^{\wedge 2}_{h\Sigma_2} \to X
	\]
is a way of encoding a coherently symmetric multiplication. In the
presence of this structure, there are natural operations
	\[
	Q^i: \mathrm{H}_*(X) \to \mathrm{H}_{*+i}X
	\]
on the mod 2 homology of $X$. These operations play a crucial
role in computations. For example, in the dual Steenrod
algebra $\mathcal{A}_*=\mathbb{F}_2[\zeta_1, \zeta_2, ...]$, 
Steinberger \cite[\S III.2]{bmms} proved that, for $i\ge 1$:
	\[
	Q^{2^i}\overline{\zeta}_i = \overline{\zeta}_{i+1}.
	\]
This computation provides key non-formal input in many results, for example:
	\begin{itemize}
	\item B\"okstedt's \cite{bokstedt}
	computation of the topological Hochschild homology
	of $\mathbb{F}_2$ and $\mathbb{Z}$.
	\item The Hopkins-Mahowald theorem identifying $\mathrm{H}\mathbb{F}_p$
	as a Thom spectrum.
	\item Lawson's \cite{lawson} proof that the Brown-Peterson spectrum
	cannot be made $\mathbb{E}_{\infty}$.
	\end{itemize}

Since the success of Hill-Hopkins-Ravenel \cite{hhr}, there
has been much recent activity in calculations in $C_2$-equivariant
homotopy theory. To aid in this endeavor, the author set out in
\cite{wilson-old} to study $C_2$-equivariant power operations on
ordinary equivariant homology with mod 2 coefficients. As first
applications, the author gave an equivariant, cellular construction of $\mathrm{BP}$
with its action by complex conjugation, and proved an analogue
of the Hopkins-Mahowald theorem jointly with Mark Behrens \cite{behrens-wilson}.
Unfortunately, the treatment
in \cite{wilson-old} is a bit clunky, contains some errors in various
formulas, and did not include the Adem relations.

The purpose of this note is to give a corrected and expanded
account of these power operations. In order to more quickly provide
a reference
for those interested in calculating, we have decided to state the main
formulas governing the operations and sketch the proofs. We defer
the verification of certain categorical statements to a more detailed,
forthcoming paper. We believe our approach is interesting even
nonequivariantly and have given a separate and complete treatment
of that story in \cite{wilson-classical}.

\subsubsection*{Outline}

In \S\ref{sec:results} we briefly recall the set-up for equivariant homotopy
theory and then go on to state the properties of the equivariant power
operations. This includes a Cartan formula (Proposition \ref{prop:cartan}),
Adem relations (Theorem \ref{thm:adem}), and Nishida relations
(Proposition \ref{prop:nishida}). We end by computing
the action of power operations on the equivariant dual Steenrod algebra
(Theorem \ref{thm:action}).

In \S\ref{sec:sketch} we sketch the
categorical and computational input necessary to prove the theorems
from \S\ref{sec:results}. The main idea is to adapt the method in
\cite{wilson-classical} to the equivariant setting using several
basic calculations combined with tools from the theory of
$C_2$-$\infty$-categories and higher equivariant algebra.

Much of this paper is devoted to explaining how the theory
of power operations in $C_2$-equivariant homotopy theory
behaves like its nonequivariant counterpart (or, more accurately,
like the \emph{odd} primary nonequivariant counterpart). 
However, there are several important instances where the theory is
different:
	\begin{itemize}
	\item The equivariant Cartan and Adem formulas contain terms which
	vanish after restriction to nonequivariant homotopy.
	\item The homology of equivariant extended powers of spheres
	is not so simple. For example, there is not, in general, a Thom
	isomorphism which reduces their computation to the homology
	of a classifying space.
	\item The power operations we build are necessarily stable.
	Nonequivariantly, all mod 2 power operations are stable. The analogous
	statement is false equivariantly, see \S\ref{ssec:counterexample}.
	\end{itemize}
	
\subsubsection*{Acknowledgements} We are grateful to Tom Bachmann,
Saul Glasman,
Jeremy Hahn, Mike Hill, Tyler Lawson, and Krishanu Sankar for 
helpful conversations. We also thank Nick Kuhn, Cary Malkiewich, Sean Tilson, 
Yan Zou, and an anonymous referee for comments on an earlier
incarnation of this material.

\subsubsection*{A historical remark} Our approach to power operations
relies heavily on their relationship to Tate constructions. This
has a long history. On the algebraic side, the relationship
between Tate-like constructions and operations goes back to Singer \cite{singer}.
Jones and Wegmann realized this algebraic construction topologically in
\cite{jones-wegmann}
(see especially their results 5.2(a) and 5.2(b)). An ad-hoc
version of the Tate-valued
Frobenius on $\mathrm{H}\mathbb{F}_p$ makes an appearance
in disguise during McClure's proof of the Nishida relations
in \cite[\S VIII.3]{bmms}, and there are hints of the same
in Steinberger's calculation of the action of power operations
on the dual Steenrod algebra \cite[\S III.2]{bmms}. More recent references,
where the relationship is more explicit,
include course notes of Lurie \cite{lurie-course}, a paper of Kuhn and McCarty
\cite[Cor. 2.13]{kuhn-mccarty}, and especially
a paper of Glasman-Lawson \cite{glasman-lawson}.
See also Theorem IV.1.15 in \cite{nikolaus-scholze} (attributed to Lurie).

\section{Results} \label{sec:results}

In \S\ref{ssec:preliminaries} we review enough equivariant
homotopy theory to be able to state our results. In
\S\ref{ssec:operations} we assert the existence of power operations
satisfying certain basic properties. In the remaining sections we
state the Cartan formulas for products and norms (Proposition
\ref{prop:cartan}), the Adem relations (Theorem \ref{thm:adem}),
the Nishida relations (Proposition \ref{prop:nishida}), and the action on the
equivariant dual Steenrod algebra (Theorem \ref{thm:action}).
In each case we give a clean
statement in terms of the total power operation and then extract
more explicit formulas from there. The formulas in terms of power
series identities are more than just a mnemonic; in \S\ref{sec:sketch}
we explain how the set-up 
we have chosen produces these identities directly. 

\begin{remark} The answers we get below closely align with 
the properties of the motivic Steenrod operations \cite{voevodsky} upon making
the substitutions:
	\[
	Q^{i\rho} \leftrightarrow \mathrm{Sq}^{-2i},
	\,\, Q^{i\rho -1} \leftrightarrow \mathrm{Sq}^{-2i+1},\,\,
	a \leftrightarrow \rho, \,\, u \leftrightarrow \tau
	\]
Of course, we allow $i$ to be positive in our case, and the proofs
are necessarily different. It would be interesting to know of a more
precise relationship between the two stories.
\end{remark}

\subsection{Preliminaries} \label{ssec:preliminaries}
We establish just as many notions and notations as we need to state our results.
A nice treatment of most of the material we need is in \cite[\S 1-3, \S A, \S B]{hhr}.
\subsubsection{Spaces and spectra}
The homotopy theory of $C_2$-spaces arises from the category
$\mathbf{Top}^{C_2}$
of topological spaces equipped with a $C_2$-action by inverting
(in the $\infty$-categorical sense)
the collection of maps $W$ which become weak equivalences on underlying and
fixed point spaces. We write:
	\[
	\mathbf{Top}^{C_2}[W^{-1}] =: \mathsf{Spaces}^{C_2}.
	\]
Given a representation $V$ we may form the one-point compactification
$S^V$ and thus define functors $\Sigma^V$ and $\Omega^V$ on
$\mathsf{Spaces}^{C_2}_*$. If we let $\sigma$ denote the sign representation
and $\rho = 1+\sigma$ the regular representation, then the
$\infty$-category of \textbf{(genuine) $C_2$-spectra} is the homotopy limit:
	\[
	\mathsf{Sp}^{C_2} = \holim\left(\cdots \to \mathsf{Spaces}^{C_2}_*
	\stackrel{\Omega^\rho}{\longrightarrow} \mathsf{Spaces}^{C_2}_*\right).
	\]
This limit comes with a canonical infinite loop functor $\Omega^{\infty}:
\mathsf{Sp}^{C_2} \to \mathsf{Spaces}^{C_2}_*$ which admits a
left adjoint, the suspension spectrum functor $\Sigma^{\infty}$.
We abuse notation and write
$S^V$ for either the space or its suspension spectrum, and
similarly if $X$ is a $C_2$-space then we often write $X_+$ for
its suspension spectrum.

The
$\infty$-category $\mathsf{Sp}^{C_2}$ is stable and hence canonically
enriched in the $\infty$-category of spectra \cite[4.2.1.33,4.8.2.18]{HA}.
In particular, the objects
$S^0$ and $C_{2+}$ corepresent functors to spectra, denoted
$(-)^{C_2}$ and $(-)^u$, respectively. We refer to $X^{C_2}$ as the
\textbf{genuine fixed points} and $X^u$ as the \textbf{underlying spectrum}.

The $\infty$-category $\mathsf{Sp}^{C_2}$ is closed
symmetric monoidal under the smash product. We denote the internal
mapping object by $F(-,-)$. 

\subsubsection{Homotopy groups}

If $V$ is a representation, then $S^V$ is an invertible object in
$\mathsf{Sp}^{C_2}$, so we can make sense of $S^{a+b\sigma}$
for any $a, b \in \mathbb{Z}$. We define groups indexed over
the representation ring $RO(C_2)$ by:
	\[
	\pi_{a+b\sigma}X := [S^{a+b\sigma}, X],
	\]
and we denote the direct sum of these groups by $\pi_{\star}X$.

By dualizing the cofiber sequence $C_{2+} \to S^0 \to S^{\sigma}$
and using the shearing equivalence $C_{2+} \wedge S^{1-\sigma} \simeq
C_{2+}$
we learn that $C_{2+}$ is self-dual. In particular, we have a map
$S^0 \to C_{2+}$ dual to the collapse $C_{2+} \to S^0$. These give
rise to the transfer and restriction, respectively, after smashing
with $S^{a+b\sigma}$:
	\[
	\mathrm{tr}_{a+b\sigma}: \pi_{a+b}X^u \to \pi_{a+b\sigma}X
	\]
	\[
	\mathrm{res}: \pi_{a+b\sigma}X \to \pi_{a+b}X
	\]
If $b=0$ then we omit the subscript on the transfer.

If $X$ is equipped with a homotopy commutative product
$X \wedge X \to X$ then the homotopy groups carry the structure
of a graded ring with commutation rule given by
\cite[Lemma 2.12]{hu-kriz}:
	\[
	xy = (-1)^{km}(1-\mathrm{tr}(1))^{\ell n}yx, \quad x\in \pi_{k+\ell\sigma}X,
	y \in \pi_{m+n\sigma}X.
	\]	
In particular, if $\mathrm{tr}(1)=2$, this is the usual rule
for a bigraded commutative ring.

\subsubsection{Eilenberg-MacLane spectra}

Associated to any abelian group $M$ there is a $C_2$-spectrum
which we will denote by the same symbol and which is uniquely characterized
by the requirements:
	\[
	\mathrm{res}: \pi_0M \stackrel{\simeq}{\to} \pi_0M^u = M
	\]
	\[
	\pi_iM=\pi_iM^u = 0, \quad i\ne0.
	\]
If the abelian group
$M$ is equipped with the structure of a commutative ring
then the corresponding $C_2$-spectrum is an $\mathbb{E}_{\infty}$-ring
in $\mathsf{Sp}^{C_2}$. We also have that $\mathrm{tr}(1) = 2$ so that
any $M$-module equipped with a homotopy commutative multiplication
satisfies the simpler bigraded commutativity rule on its homotopy groups.

\begin{convention} We will denote by $k$ both the field $\mathbb{F}_2$
and its corresponding $C_2$-spectrum.
\end{convention}

We will need two elements in the homotopy of $k$. The first comes from
the homotopy of the sphere:
	\[
	a: S^{-\sigma} \to S^0
	\]
is adjoint to the inclusion of fixed points. The second arises from the identity
$\mathrm{tr}(1) = 0 \in \pi_0k$, which produces a map 
	\[
	u: S^{1-\sigma} =
	\mathrm{cofib}(S^0 \to C_{2+}) \to k
	\]
that is well-defined since $\pi_1k = 0$. 

\subsubsection{Norms}

The symmetric monoidal structure on $\mathsf{Sp}^{C_2}$
can be enhanced to allow for tensor products indexed over finite $C_2$-sets.
This is one of the key innovations introduced in \cite[\S B]{hhr}.

To describe the homotopical features of this enhancement, we define,
for each $n$, an $\infty$-category $\left(\mathsf{Sp}^{C_2}\right)^{h\Sigma_{[n]}}$
as follows. Begin with $\mathbf{Top}^{C_2 \times \Sigma_n}$ and invert
those maps which are $C_2$-equivalences to form an $\infty$-category
$\left(\mathsf{Spaces}^{C_2}\right)^{h\Sigma_{[n]}}$. Now form the limit
	\[
	\left(\mathsf{Sp}^{C_2}\right)^{h\Sigma_{[n]}}:=
	\holim \left(\cdots \stackrel{\Omega^{\rho}}{\longrightarrow}
	\left(\mathsf{Spaces}^{C_2}\right)_*^{h\Sigma_{[n]}}
	\stackrel{\Omega^{\rho}}{\longrightarrow}
	\left(\mathsf{Spaces}^{C_2}\right)_*^{h\Sigma_{[n]}}\right).
	\]
The inclusion $\mathbf{Top}^{C_2} \to \mathbf{Top}^{C_2\times\Sigma_n}$
of $C_2$-spaces with trivial $\Sigma_n$-action induces functors
	\[
	\mathsf{Spaces}^{C_2} \longrightarrow
	\left(\mathsf{Spaces}^{C_2}\right)^{h\Sigma_{[n]}}
	\]
	\[
	\mathsf{Sp}^{C_2} \longrightarrow
	\left(\mathsf{Sp}^{C_2}\right)^{h\Sigma_{[n]}}
	\]
which admit both left and right adjoints, denoted $(-)_{h\Sigma_{[n]}}$
and $(-)^{h\Sigma_{[n]}}$, respectively.

We will need to know that the assignment $X \mapsto X^{\wedge n}$ can be
refined to a functor
	\[
	(-)^{\wedge n}: \mathsf{Sp}^{C_2}\longrightarrow
	\left(\mathsf{Sp}^{C_2}\right)^{h\Sigma_{[n]}},
	\]
which allows us to produce a homotopical symmetric power
	\[
	\mathrm{Sym}^{[n]}(X):= (X^{\wedge n})_{h\Sigma_{[n]}}.
	\]
The latter functor, at least, can be found in \cite[\S B.6.1]{hhr}.
In the case $n=2$ we may form the \textbf{norm}:
	\[
	N(X):= \left(X^{\wedge 2} \wedge 
	\frac{C_2\times \Sigma_2}{\Delta}_+\right)_{h\Sigma_{[2]}} \in 
	\mathsf{Sp}^{C_2}.
	\]
Note that we have natural maps of $C_2$-spectra:
	\[
	X^{\wedge 2} \to \mathrm{Sym}^{[2]}(X), \quad
	N(X) \to \mathrm{Sym}^{[2]}(X).
	\]

We will also need the generalization of the above discussion to
$k$-modules, where the only difference is that we use
$\otimes := \otimes_k$ instead of the smash product.

\subsection{The operations} \label{ssec:operations}

The most general setting in which our operations exist
is the following.

\begin{definition} A $k$-module $A$ is \textbf{equipped with
an equivariant symmetric multiplication} if we provide a map
$\mathrm{Sym}^{[2]}(A) \to A$, where we take the indexed
symmetric square in $\mathsf{Mod}_k$.
\end{definition}

Note that we do not require this product to be unital or associative.
An example of such a $k$-module is $F(X_+, k)$ where $X$ is a
$C_2$-space and the product arises from the diagonal.

\begin{definition} If $A_{\star}$ is an $RO(C_2)$-graded commutative
$k_{\star}$-algebra, we denote by $A_{\star}(\!(s,t)\!)$ the algebra:
	\[
	A_{\star}(\!(s,t)\!):= (A_{\star}[s]\llbracket t\rrbracket/(s^2=as+ut))[t^{-1}]
	\]
where $|s|=-\sigma$ and $|t|=\rho$.
We observe that homogeneous
elements in this algebra may be written uniquely in the form
	\[
	\sum_{i\in \mathbb{Z}} x_i st^{i-1} + \sum_{i\in \mathbb{Z}} y_i t^{i}
	\]
where $x_i, y_i \in A_{\star}$ vanish for $i$ sufficiently negative.
\end{definition}

In Construction \ref{cstr:operations} below, we will sketch the construction
of operations $Q^{i\rho -\varepsilon}: A_{\star} \to A_{\star+i\rho-\varepsilon}$
for $\varepsilon = 0,1$. These agree with the operations constructed
in a different manner in our previous treatment \cite{wilson-old}, and
the following theorem summarizes their first properties. We draw attention
to (iv)-(vi) which are corrected versions of the corresponding statements
in \cite{wilson-old}.

\begin{theorem}\label{thm:properties} Let $A$ be a $k$-module equipped with
an equivariant symmetric multiplication. Then there is
a natural operation:
	\[
	Q(s,t): A_{\star} \to A_{\star}(\!(s,t)\!)
	\]
with coefficients
	\[
	Q(s,t)x = \sum_i (Q^{i\rho-1}x) st^{i-1} + \sum_i (Q^{i\rho}x) t^{i}
	\]
satisfying the following properties:
	\begin{enumerate}[\upshape(i)]
	\item \textup{(Mackey)} For $\varepsilon=0,1$, the operation
	$Q^{i\rho-\varepsilon}$ restricts to $Q^{2i-\varepsilon}$ on
	underlying homotopy groups, and commutes with addition and transfer.
	\item \textup{(Loops)} For any representation $V$, we have
	$\Omega^V Q(s,t) = Q(s,t) \Omega^V$.
	\item \textup{(Squaring)} If $|x| = n\rho-\varepsilon$ for $\varepsilon=
	0,1$ then $Q^{n\rho-\varepsilon}(x) = x^2$.
	\item \label{vanishing} \textup{(Vanishing)} Suppose $|x| = a+b\sigma$
	and $\varepsilon = 0,1$. Then
	$Q^{i\rho-\varepsilon}x = 0$ if $i<a+\varepsilon$ and $i\le b$.
	\item \textup{(Cohomology)} If $A=F(X,k)$ for a pointed $C_2$-space
	$X$,
	then $Q^{i\rho-\varepsilon}$ acts by zero if $2i-\varepsilon>0$.
	\item \textup{(Bockstein)} Suppose that $A$ arises as the
	mod 2 reduction of a $\mathbb{Z}/4$-module. Then
	$Q^{i\rho -1}$ acts by $\beta Q^{i\rho}$ where $\beta$ is the Bockstein.
	\end{enumerate}
\end{theorem}

\subsection{Cartan formulas} \label{ssec:cartan}

In \cite{wilson-old} we stated an obviously incorrect version
of the Cartan formula (it did not even restrict to the usual
version on underlying homotopy). The correct version
of the Cartan formula, as well as a new formula for the 
value of operations on a norm class, is as follows:

\begin{proposition}\label{prop:cartan} 
	\begin{align*}
	Q(s,t)(x\otimes y) &= Q(s,t)x \otimes Q(s,t)y\\
	Q(s,t)(Nx) &= N(Q(w)x)
	\end{align*}
\end{proposition}

Here we have used $Nx$ to denote the class obtained from $x: \Sigma^ik^u
\to A^u$ as $\Sigma^{i\rho}k = N(\Sigma^ik^u) \to N(A^u)$.
We also use $w$ to denote a power series generator
for for the underlying total power operation $Q(w) = \sum Q^i(-) w^i$ and
define $N(w) = t$. 

We may extract more explicit formulas by comparing coefficients
of $s^{\varepsilon}t^{n-\varepsilon}$ and using the distributive
law for the norm of a sum \cite[Prop. A.37]{hhr}:

\begin{corollary}[Explicit Cartan formulas]
	\begin{align*}
	Q^{n\rho}(x \otimes y) &=
	\sum_{i+j=n} Q^{i\rho}x \otimes Q^{j\rho}y
	+ u\sum_{i+j=n+1}Q^{i\rho-1}x \otimes Q^{j\rho-1}y\\
	Q^{n\rho-1}(x \otimes y) &=
	\sum_{i+j = n} \left(Q^{i\rho -1}x \otimes Q^{j\rho}y +
	Q^{i\rho}x \otimes Q^{j\rho -1}y\right) +
	a\sum_{i+j = s+1} Q^{i\rho-1}x \otimes Q^{j\rho-1}y\\
	Q^{n\rho}(Nx) &= N(Q^nx) + \sum_{i<j, i+j = 2n}
	\mathrm{tr}_{n\rho}(Q^ix \otimes Q^jx)\\
	Q^{n\rho -1}(Nx) &=\sum_{i<j, i+j=2n-1} \mathrm{tr}_{n\rho-1}
	(Q^ix\otimes Q^jy)
	\end{align*}
\end{corollary}

There are internal versions of these formulas if the maps
$N(A^u) \to \mathrm{Sym}^{[2]}(A) \to A$ and
$A^{\otimes 2} \to \mathrm{Sym}^{[2]}(A) \to A$ commute
with the equivariant symmetric multiplication.

\subsection{Adem relations} \label{ssec:adem}

The author's motivation for the change in our
treatment of power operations was to prove the Adem relations. We state
them first in a form reminscent of Bullett-Macdonald \cite{bullett-macdonald},
Steiner \cite{steiner}, and Bisson-Joyal \cite{bisson-joyal},
but we need some preliminaries.

\begin{definition} An \textbf{Adem object} is a $k$-module
$A$ equipped with an equivariant symmetric multiplication
such that there exists a dotted arrow
	\[
	\xymatrix{
	\mathrm{Sym}^{[2]}(\mathrm{Sym}^{[2]}(A)) \ar[r]\ar[d]&
	\mathrm{Sym}^{[2]}(A) \ar[d]\\
	\mathrm{Sym}^{[4]}(A) \ar@{-->}[r] & A
	}
	\]
making the diagram commute up to homotopy.
\end{definition}

\begin{definition} Given $Q(s,t): A_{\star} \to
A_{\star}(\!(s,t)\!)$ we extend to a ring homomorphism
	\[
	Q(s,t): A_{\star}(\!(c,d)\!) \to A_{\star}(\!(c,s,d,t)\!)
	\]
by the rule:
	\[
	Q(s,t)c = c+dst^{-1}, \quad Q(s,t)d = d+d^2t^{-1}.
	\]
Here the target ring is defined as:
	\[
	A_{\star}(\!(c,s,d,t)\!) :=
	\left(A_{\star}[c,d]\llbracket d,t\rrbracket/(c^2=ac+ud,
	s^2=as+ut)\right)[d^{-1},t^{-1}]
	\]
and the operation lands in this ring because of the vanishing
property in Theorem \ref{thm:properties}(\ref{vanishing}).
\end{definition}

\begin{theorem}[Adem relations] \label{thm:adem} Suppose $A$ is an Adem object.
Then, for any $x \in A$, the series
	\[
	Q(s,t) \circ Q(c,d) x \in A_{\star}(\!(c,s,d,t)\!)
	\]
is invariant under the transformation $c \leftrightarrow s$,
$d \leftrightarrow t$.
\end{theorem}

Comparing coefficients of $s^{\varepsilon}c^{\varepsilon'}$ we get
four different identities between series in $d$ and $t$. Using
a residue argument analogous to that in \cite{bullett-macdonald}
or \cite[\S 1]{bisson-joyal}, we can extract the following:

\begin{corollary}[Explicit Adem relations] With assumptions as above,
we have (with $x$ omitted):
	\begin{align*}
	Q^{i\rho -1}Q^{j\rho -1} &=
	\sum \binom{\ell-j-1}{2\ell-i}Q^{(i+j-\ell)\rho-1}Q^{\ell \rho -1}
	\\
	Q^{i\rho}Q^{j\rho} &=
	\sum \binom{\ell-j-1}{2\ell-i}Q^{(i+j-\ell)\rho} Q^{\ell\rho}
	+ u \sum 
	\binom{\ell-j-1}{2\ell-1-i}Q^{(i+j-\ell+1)\rho -1}Q^{\ell\rho -1}
	\\
	Q^{i\rho -1}Q^{j\rho}&=
	\sum \binom{\ell-j-1}{2\ell-i}Q^{(i+j-\ell)\rho -1}Q^{\ell\rho} +
	a\sum \binom{\ell-j-1}{2\ell-i-1}Q^{(i+j-\ell+1)\rho-1}Q^{\ell\rho-1}\\
	Q^{i\rho}Q^{j\rho -1} &=
	\sum_{\ell} \binom{\ell-j}{2\ell-i} Q^{(i+j-\ell)\rho-1}Q^{\ell\rho}
	+\sum_{\ell}\binom{\ell-j}{2\ell+1-i}Q^{(i+j-\ell-1)\rho}Q^{(\ell+1)\rho-1}
	\\&+ a\sum \binom{\ell-j}{2\ell+1-i} Q^{(i+j-\ell)\rho-1}Q^{(\ell+1)\rho-1}
	\end{align*}
\end{corollary}

\subsection{Nishida relations} \label{ssec:nishida}

Recall that Hu and Kriz \cite{hu-kriz} computed the equivariant
dual Steenrod algebra to be
	\[
	\mathcal{A}^{C_2}_{\star} := \pi_{\star}(k \wedge k)
	= k_{\star}[\tau_i, \xi_{i+1}| i\ge 0]/(\tau_i^2 = (u+a\tau_0)\xi_{i+1} +
	a\tau_{i+1})
	\]
for certain generators with
	\[
	|\tau_i| = 2^i\rho -\sigma, \quad |\xi_i| = (2^i-1)\rho.
	\]
If $X$ is a $C_2$-spectrum then $k \wedge X$ is a equipped with a
right coaction
	\[
	k_{\star}X \to (k_{\star}X) \otimes_{k_{\star}} \mathcal{A}^{C_2}_{\star}
	\]
arising from the unit $k \wedge X \to k \wedge X \wedge k$. The Nishida
relations describe the relationship between power operations on
$k_{\star}X$, if they exist, and this coaction.

The Nishida relations from \cite{wilson-old} are correct as stated,
but we restate them here in a slightly different form. We will use formal
series:
	\[
	\tau(s,t) = s + \sum_{i\ge 0} \tau_i t^{2^i},
	\,\, \xi(t) = t+ \sum_{i\ge 1}t^{2^i}
	\]
as well as similar series for the conjugates 
$\overline{\tau}(s,t)$ and $\overline{\xi}(t)$. Note
$\tau(s,t)+s$ is a power series in $t$.

\begin{proposition}\label{prop:nishida} Suppose $X$ is a $C_2$-spectrum
equipped with an equivariant symmetric multiplication
and $x \in \pi_{\star}(k \wedge X)$ is arbitrary. Then
	\[
	\sum \psi_R(Q^{i\rho-1}x) st^{i-1} +
	\sum \psi_R(Q^{i\rho}x) t^{i} =
	Q(\overline{\tau}(s,t), \overline{\xi}(t))\psi_R(x).
	\]
\end{proposition}

\begin{remark} The theorem applies more generally to
limits of spectra of the form $k \wedge X$ equipped with
the corresponding \emph{completed} right coaction,
e.g. to $k^{h\Sigma_{[2]}}$.
\end{remark}

Comparing coefficients yields the formulas:

\begin{corollary} With assumptions as above,
	\[
	\sum_i \psi_R(Q^{i\rho-1}x) t^{i-1} = \sum_r Q^{r\rho-1}(\psi_R(x))
	\overline{\xi}(t)^{r-1}
	\]
	\[
	\sum_i \psi_R(Q^{i\rho}x) t^i =
	\sum_r Q^{r\rho}(\psi_r(x))\overline{\xi}(t)^r +
	\sum_r Q^{r\rho-1}(\psi_R(x)) (\overline{\tau}(s,t) + s)\overline{\xi}(t)^{r-1}
	\]
\end{corollary}

\subsection{Action on dual Steenrod algebra} \label{ssec:action}

In \cite{wilson-old} we gave a partial computation of the 
action of the power operations on the dual Steenrod algebra.
The proof given was incomplete, but the results are correct.
Now we are able to state a \emph{complete} calculation of the
action.

\begin{theorem}\label{thm:action} The total power operation on 
$\mathcal{A}^{C_2}_{\star}$
is uniquely determined by the identities:
	\[
	\tau(c,d) + \xi(d)\tau(s,t)\xi(t)^{-1} = (c+dst^{-1})+\sum_{i\ge 0}
	(Q(s,t)\tau_i)(d^{2^i}+d^{2^{i+1}}t^{-2^i})
	\]
	\[
	\xi(d) + \xi(d)^2\xi(t)^{-1} =
	\sum_{i\ge 0}(Q(s,t)\xi_i)(d^{2^i} + d^{2^{i+1}}t^{-2^i})
	\]
\end{theorem}

Comparing coefficients leads to a recursive description of 
$Q(s,t)\tau_i$ and $Q(s,t)\xi_i$ which can be solved explicitly.
We record this solution and a few other more explicit corollaries
of the above formula, including the main result from \cite{wilson-old}.

\begin{corollary} The action of power operations on
$\tau_n$ and $\xi_n$ is given by the identites:
	\begin{align*}
	t^{2^n}\sum_{r} Q^{r\rho}(\tau_n) t^r &=
	\left(\sum_{i\ge n+1}\tau_i t^{2^i}\right) +
	(\tau(s,t)+s)\xi(t)^{-1}\left(\sum_{i\ge n+1}\xi_i t^{2^i}\right)\\
	t^{2^n}\sum_r Q^{r\rho -1}(\tau_n) t^{r-1} &=
	\xi(t)^{-1}\left(\sum_{i\ge n+1}\xi_i t^{2^i}\right)\\
	t^{2^n}\sum_r Q^{r\rho}(\xi_n)t^r &= \sum_{i\ge n+1}\xi_i t^{2^i}
	+\xi(t)^{-1}\sum_{i\ge n} \xi_i^2t^{2^{i+1}}\\
	\sum_r Q^{r\rho-1}(\xi_n)t^{r-1} &= 0
	\end{align*}
\end{corollary}

\begin{corollary} 	
	The following formulas, for $k\ge 1$, hold for the action on
	$\tau_0$:
	\begin{align*}
	Q^{(2^k-1)\rho}\tau_0 &= \overline{\tau}_k, \\
	Q^{(2^k-1)\rho -1}\tau_0 &= \overline{\xi}_k.
	\end{align*}
\end{corollary}

\begin{corollary}
	The following formulas hold for the action on
	$\tau_i, \xi_i$, and their conjugates:
	\begin{align*}
	Q^{2^k\rho}\tau_k &= \tau_{k+1} + \tau_0\xi_{k+1}
	& Q^{2^k\rho}\overline{\tau}_k &=\overline{\tau}_{k+1} \\
	Q^{2^k\rho -1}\tau_k &= \xi_{k+1}
	& Q^{2^k\rho -1}\overline{\tau}_k &=\overline{\xi}_{k+1}\\
	Q^{2^k\rho}\xi_k &= \xi_{k+1} + \xi_1\xi_{k}^2
	& Q^{2^k\rho}\overline{\xi}_k &= \overline{\xi}_{k+1}
	\end{align*}
\end{corollary}

As an algebraic corollary, used in \cite{behrens-wilson}, we have:

\begin{corollary} The algebra $\mathcal{A}^{C_2}_{\star}$ is generated
by $\tau_0$ as a ring equipped with the operations $Q^{i\rho-\varepsilon}$.
More specifically, the elements $Q^{2^k\rho}Q^{2^{k-1}\rho}\cdots Q^{\rho}\tau_0$
and their Bocksteins generate $\mathcal{A}^{C_2}_{\star}$.
\end{corollary}

\section{Sketch of the proof}\label{sec:sketch}

\subsection{Review of the classical case}\label{ssec:review}
We begin by reviewing the classical setting, with some modifications.
It is possible \cite{wilson-classical}
to develop the entire theory of mod 2 power operations based
on the following categorical inputs:
	\begin{enumerate}[(A)]
	\item (Universal property of the Tate construction) 
	\cite[\S I.3]{nikolaus-scholze}
	Let $\mathcal{F}$ be a collection of subgroups of
	a group $G$, closed under sub-conjugacy. Then the
	generalized Tate construction \cite[p.443]{greenlees1}
		\[
		(-)^{t\mathcal{F}}: \mathsf{Sp}^{hG} \to \mathsf{Sp}
		\]
	has the following universal property: the map
	$(-)^{hG} \to (-)^{t\mathcal{F}}$ is initial amongst
	all natural transformations of exact, lax symmetric monoidal
	functors whose target annihilates the stable subcategory
	of objects induced from subgroups 
	$H \in \mathcal{F}$.
	\item (Universal property of the forgetful functor)
	\cite{nikolaus}
	The forgetful functor $U: \mathsf{Mod}_{\mathbb{F}_2}
	\to \mathsf{Sp}$ is initial amongst all exact, lax symmetric
	monoidal functors from $\mathsf{Mod}_{\mathbb{F}_2}$
	to $\mathsf{Sp}$.
	\item (Universal property of Tate powers) Let $G \subseteq \Sigma_n$
	be a subgroup and $\mathcal{T}$ the family of non-transitive subgroups
	of $G$. Define the Tate powers by
	$T_G(X):= (X^{\otimes n})^{t\mathcal{T}}$ and the divided power by
	$\Gamma^G(X):= (X^{\otimes n})^{hG}$. 
	Then the map $\Gamma^G \to T_G$ is initial amongst natural transformations
	from $\Gamma^G$ to an exact functor, i.e. $T_G$ is the Goodwillie derivative
	of $\Gamma^G$ and is, in particular, exact. We do not
	know a reference that states the result in these words, but the
	actual calculation is well-known; we spell this out in
	\cite[2.2.1]{wilson-classical}.
	\end{enumerate}
We also need the following computational inputs:
	\begin{enumerate}[(I)]
	\item As a ring, 
		\[
		\pi_*\mathbb{F}_2^{t\Sigma_2} = \mathbb{F}_2(\!(w)\!)
		\]
	where $w$ is the first Stiefel-Whitney class.
	\item There is a canonical
	equivalence $(\Sigma^{k}\mathbb{F}_2)^{\otimes 2}
	\simeq \Sigma^{2k}\mathbb{F}_2$ in $\mathsf{Mod}_{\mathbb{F}_2}^{h\Sigma_2}$.
	\item The map $\mathbb{F}_2 = \mathbb{F}_2 \wedge S^0 \to
	\mathbb{F}_2 \wedge \mathbb{F}_2$ induces, upon taking Tate
	fixed points for the trivial action of $\mathbb{F}_2$, the completed
	coaction:
		\[
		\psi_R: \mathbb{F}_2(\!(w)\!) \to \mathcal{A}_*(\!(w)\!)
		\]
	determined by 
		\[
		\psi_R(w) = \sum \zeta_i w^{2^i} =: \zeta(w).
		\]
	\end{enumerate}
We now review how these facts imply the desired properties of
power operations. A more detailed treatment is given in \cite{wilson-classical}.
	
Suppose we are given an $\mathbb{F}_2$-module $A$ equipped
with a symmetric multiplication $\mathrm{Sym}^2(A) \to A$. By
(A), (B), and (C) there is an essentially unique transformation of
lax symmetric monoidal functors $U \to UT_2$, and so we may form
the \textbf{Tate-valued Frobenius}:
	\[
	A \to T_2(A) = (A^{\otimes 2})^{t\Sigma_2} \to A^{t\Sigma_2}.
	\]
We can then define operations $Q^i: A \to \Sigma^{-i}A$,
using (I), as the composite
	\[
	A \to A^{t\Sigma_2} \stackrel{w^{-i-1}}{\longrightarrow}
	\Sigma^{-i-1}A^{t\Sigma_2} \to \Sigma^{-i}A_{h\Sigma_2} \to \Sigma^{-i}A
	\]
The Frobenius then induces a total power operation
	\[
	Q(w): A_* \to A^{t\Sigma_2}_* \simeq A_*(\!(w)\!)
	\]
given as $Q(w)x = \sum_i (Q^ix) w^i$. The first properties follow
immediately from this definition, with the exception of the squaring
property, which eventually requires the computation (II) above.
The Cartan formula follows from the lax symmetric monoidal structures
on every functor in sight. The Adem relations follow from an important
lemma: every $\mathbb{F}_2$-module equipped with a symmetric multiplication
admits a canonical lift of the iterated Frobenius
	\[
	\xymatrix{
	&A^{\tau\Sigma_2\wr \Sigma_2}\ar[d]\\
	A \ar[r]\ar[ur] & (A^{t\Sigma_2})^{t\Sigma_2}
	}
	\]
where $(-)^{\tau G}$ denotes $(-)^{t\mathcal{T}}$ when $G \subseteq \Sigma_n$
is a subgroup. The proof of this uses the universality assertions,
naturality, and lax symmetric monoidal properties in each of
(A), (B), and (C) to verify that various diagrams commute.

If $A$ is, moreover, an Adem object, then there
is a further lift to $A^{\tau\Sigma_4}$ and this immediately implies
the Adem relations in the form stated by \cite{steiner} after Bullett-Macdonald
\cite{bullett-macdonald}. 

The computation (III) immediately implies the Nishida relations
by examining the diagram
	\[
	\xymatrix{
	(\mathbb{F}_2 \wedge X) \ar[r] \ar[d] & (\mathbb{F}_2 \wedge X)
	\otimes (\mathbb{F}_2 \wedge \mathbb{F}_2) \ar[d]\\
	(\mathbb{F}_2 \wedge X)^{t\Sigma_2} \ar[r] &
	\left((\mathbb{F}_2 \wedge X)
	\otimes (\mathbb{F}_2 \wedge \mathbb{F}_2)\right)^{t\Sigma_2}
	}
	\]
	
Finally, the action on the dual Steenrod algebra is obtained by
applying the Nishida relations, in the form above, to the 
element $w \in \pi_{-1}\mathbb{F}_2^{h\Sigma_2}$. 

\subsection{The equivariant setting: categorical input}

The $C_2$-equivariant setting proceeds in an identical fashion.
The idea is to first make the following adjustments to the
categorical input:
	\begin{itemize}
	\item Replace the use of $\infty$-categories with
	$C_2$-$\infty$-categories, in the sense of
	\cite{bdgns}.
	\item Replace the use of exact functors with
	$C_2$-exact functors in the sense of \cite[2.14]{nardin}.
	\item Replace the use of symmetric monoidal
	$\infty$-categories and lax symmetric monoidal functors
	with their indexed analogues. One can define
	a $C_2$-symmetric monoidal $\infty$-category
	as a product preserving functor
		\[
		\mathsf{Span}(\mathsf{Fin}^G) \to \mathsf{Cat}_{\infty}
		\]
	but the notion of lax $C_2$-symmetric monoidal functors
	takes a little longer to define.
	A suitable definition can be found in
	\cite[2.1.4.2]{wilson-thesis}. 
	\item Replace the group $\Sigma_n$ with the
	$C_2$-groupoid $\Sigma_{[n]}$, i.e. the functor
	$\Sigma_{[n]}: \mathrm{Orbit}(C_2)^{\mathrm{op}} \to
	\mathsf{Gpd}$ which assigns to an orbit $T$ the
	groupoid of maps of $(C_2\times \Sigma_n)$-sets
	$U \to T$ which are $\Sigma_n$ torsors. This records
	the natural symmetries on indexed tensor powers
	in a $C_2$-symmetric monoidal $\infty$-category.
	\end{itemize}
With these changes, the claims (A), (B), and (C) remain true
in the $C_2$-equivariant setting, and the proofs are mostly
the same as their nonequivariant analogues. We note that,
without the symmetric monoidal structures, these statements
are readily provable using tools from \cite{nardin} and \cite{shah}. 
However, the additional multiplicative structure requires more set-up.

\subsection{The equivariant setting: computational input}

The computation (I) is replaced
by 
	\[
	\pi_{\star}k^{t\Sigma_{[2]}} = k_{\star}(\!(s,t)\!)
	\]
with notation as in \S\ref{ssec:operations}. Here some care is required.
First, recall that there is a $C_2$-space $\mathrm{B}_{C_2}\Sigma_2$
which classifies $\Sigma_2$-torsors on $C_2$-spaces,
and a corresponding universal torsor $\mathrm{E}_{C_2}\Sigma_2$.
The map $\mathrm{E}_{C_2}\Sigma_{2+} \to S^0$ is an equivalence
in $\left(\mathsf{Sp}^{C_2}\right)^{h\Sigma_{[2]}}$, so we get
equivalences:
	\[
	F(\mathrm{E}_{C_2}\Sigma_{2+}, X)^{h\Sigma_{[2]}} \simeq
	X^{h\Sigma_{[2]}}
	\]
	\[
	(X \wedge \mathrm{E}_{C_2}\Sigma_{2+})_{h\Sigma_{[2]}}
	\simeq X_{h\Sigma_{[2]}}
	\]
for any $X \in \left(\mathsf{Sp}^{C_2}\right)^{h\Sigma_{[2]}}$. We will
use this to define a filtration on these functors.

If we let
$\tau$ denote the sign representation of $\Sigma_2$, then an explicit
model for $\mathrm{E}_{C_2}\Sigma_2$
is the unit sphere $S(\infty(\rho\otimes \tau))$
in the $C_2\times\Sigma_2$-representation $(\rho\otimes \tau)^{\oplus \infty}$.
The filtration
	\[
	\Sigma_{2+} = S(1\otimes \tau) \subseteq
	S(\rho \otimes \tau) \subseteq S((\rho+1)\otimes \tau)
	\subseteq \cdots 
	\]
then produces a filtration on the functors $(-)^{h\Sigma_{[2]}}$ and
$(-)_{h\Sigma_{[2]}}$. Now, in general, if $\lambda$ is a one-dimensional
representation and $W$ is any representation, then the map
	\[
	S(\lambda) \times \mathrm{Hom}(\lambda, W) \to
	S(\lambda \oplus W) - S(W), \quad (x,\phi) \mapsto
	\frac{(x,\phi(x))}{|(x,\phi(x))|}
	\]
is an equivariant homeomorphism and so induces an equivalence
	\[
	S(\lambda)_+ \wedge S^{\lambda^{\vee}\otimes W} \simeq
	S(\lambda \oplus W)/S(W)
	\]
upon taking one-point compactifications. Applying this to $S(\tau)=\Sigma_2$
and $S(\sigma\otimes \tau) = \frac{C_2\times\Sigma_2}{\Delta}$, where
$\Delta$ is the diagonal subgroup, we learn that:
	\[
	\mathrm{gr}_{2j}\mathrm{Sym}^{[2]}(M) \simeq \Sigma^{j\rho}M^{\otimes 2},
	\]
	\[
	\mathrm{gr}_{2j+1}\mathrm{Sym}^{[2]}(M) \simeq \Sigma^{(j+1)\rho -1}
	N(M).
	\]
Unlike the classical setting, this filtration does not split in general.
For example, it does not split when $M$ is of the form $\Sigma^{j\rho +1}k$
or $\Sigma^{j\rho -2}k$ (see below \S\ref{ssec:counterexample}).

However, in the case $M=k$ or, more generally, $M=\Sigma^{j\rho-\varepsilon}k$
for $\varepsilon=0,1$, the filtration splits for degree reasons stemming
from a computation of the homology of a point, $k_{\star}$.
This gives an additive calculation
	\[
	k_{h\Sigma_{[2]}} \simeq \bigoplus_{j\ge 0}
	\Sigma^{j\rho}k \oplus \bigoplus_{j\ge 1}
	\Sigma^{j\rho-1}k
	\]
Dualizing produces an additive calculation for $k^{h\Sigma_{[2]}}$ and it
is not difficult to prove that the ring structure is given by
	\[
	k^{h\Sigma_2}_{\star}= k_{\star}[s]\llbracket t\rrbracket/(s^2=as+ut)
	\]
This calculation was obtained in a different way by Hu-Kriz \cite[6.27]{hu-kriz},
but the filtration above is useful in other contexts.

\begin{remark} Here there is a subtle point: there are two possible generators in degree
$-\sigma$ which satisfy the above relation.
The choice is determined by restriction to
$S(\rho\otimes \tau)/\Sigma_2 = S^{\sigma} \subseteq \mathrm{B}_{C_2}\Sigma_2$.
It is standard to choose the generator which restricts to 
$1 \in \pi_0k =\pi_{-\sigma}\Sigma^{-\sigma}k$. We will \emph{not} make that choice.
Instead we choose $s$ to be the generator which restricts to $a$. The other
generator is given by
	\[
	\overline{s}:= s+a.
	\]
These two generators behave in an almost identical fashion, but with our
choice of generator the map
	\[
	k^{t\Sigma_{[2]}} \to \Sigma k_{h\Sigma_{[2]}} \to \Sigma k
	\]
will, on homotopy, record the coefficient of $st^{-1}$ (as opposed to $(s+a)t^{-1}$),
which is more convenient.
\end{remark}

To finish the analogue of (I), one checks that the Tate construction inverts the
class $t$. The filtration can also be used to show that, for $M$ equipped with
trivial $\Sigma_2$-action,
	\[
	M^{t\Sigma_{[2]}}_{\star} = M_{\star}(\!(s,t)\!)
	\]
even if $M$ is not free as a $k$-module. An equivariant symmetric multiplication
on a $k$-module $A$ is now enough to form the
\textbf{equivariant Tate-valued Frobenius}:
	\[
	A \to T_{[2]}(A):= (A^{\otimes 2})^{t\Sigma_{[2]}} \to A^{t\Sigma_{[2]}}.
	\]
\begin{construction}\label{cstr:operations}
The individual power operations are defined as the composites:
	\[
	Q^{i\rho -\varepsilon}:
	A \to A^{t\Sigma_{[2]}} \stackrel{\overline{s}^{\varepsilon}t^{-i-1+\varepsilon}}{\longrightarrow}
	\Sigma^{-i\rho -1+\varepsilon}A^{t\Sigma_{[2]}} \to \Sigma^{-i\rho+\varepsilon}
	A_{h\Sigma_{[2]}} \to \Sigma^{-i\rho+\varepsilon}A.
	\]
\end{construction}

The analogous computation to (II) is an equivalence
$(\Sigma^{j\rho}k)^{\otimes 2} \simeq \Sigma^{2j\rho}k$ in
$\mathsf{Mod}^{h\Sigma_{[2]}}_k$. One can show this either
by establishing the orientability of a bundle on $\mathrm{B}_{C_2}\Sigma_2$
or using the computation of $\mathrm{Sym}^2(\Sigma^{j\rho}k)$ referenced above.
The analogous fact does not hold for other representations, but this
turns out to suffice. Finally, the analogue of (III) is essentially the definition
of the Milnor generators from Hu-Kriz \cite[p.381]{hu-kriz}. The main theorems
in \S\ref{sec:results} now follow as in \S\ref{ssec:review}.

\subsection{A counterexample}\label{ssec:counterexample}

In the previous section we constructed a filtration on $\mathrm{Sym}^{[2]}(M)$
and remarked that it does not split in general. Let us explain what happens
in two examples to see how things differ from the classical case.

First, in the case $M=\Sigma^{n\rho +1}k$ it turns out that, while the
given filtration does not split, we can still describe the answer:
\begin{proposition} 
	\[
	\mathrm{Sym}^{[2]}(\Sigma^{n\rho+1}k)
	\simeq
	\bigoplus_{j\ge n+1} \Sigma^{(n+j)\rho+1}k \oplus
	\bigoplus_{j\ge n+1}\Sigma^{(n+j)\rho}k
	\]
\end{proposition}

This situation is not too bad because the generators can actually
be described in terms of the power operations we have.
That is, we can choose generators so that, if $A$ is equipped with
an equivariant symmetric multiplication, and $x: \Sigma^{n\rho+1}k \to A$
is a class, then the image of the generator in degree $(n\rho+1)+(j\rho -\varepsilon)$
under the map
	\[
	\mathrm{Sym}^{[2]}(\Sigma^{n\rho+1}k) \to \mathrm{Sym}^{[2]}(A) \to A
	\]
is equal to $Q^{j\rho-\varepsilon}x$.

We also have a computation of $\mathrm{Sym}^{[2]}(\Sigma^{n\rho -2}k)$:
\begin{proposition} 
	\[
	\mathrm{Sym}^{[2]}(\Sigma^{n\rho-2}k)
	\simeq \Sigma^{2n\rho-4}k \oplus 
	\bigoplus_{j\ge n} \Sigma^{(n+j)\rho -2}k \oplus
	\bigoplus_{j\ge n}\Sigma^{(n+j)\rho-3}k
	\]
\end{proposition}

In this case, most classes will correspond to named power operations,
but the bottom class cannot, for degree reasons. In other words:
there is no stable operation which refines the square $x \mapsto x^2$
for classes in degree $n\rho -2$. We imagine the situation becomes
worse the further one strays from regular representations.

\bibliographystyle{amsalpha}
\bibliography{Bibliography}

\end{document}